\documentclass[english,11pt,a4paper]{article}
\usepackage{babel,indentfirst}
\begin{document}

\title{Solution of the quasipositivity problem in the braid group}
\author{ A.Bentalha}

\maketitle

\paragraph{Introduction}

Let $B_{n+1}$ be the group of braids with $n+1$ strings.
It has a presentation  ,
$$(\Theta): B_{n+1}= [{\sigma_1,\dots,\sigma_n\,|\, 	\sigma_j\sigma_k=\sigma_k\sigma_j\  for  |j-k|>1,
\sigma_j\sigma_k\sigma_j=\sigma_k\sigma_j\sigma_k\  for  |j-k|=1}]$$

A braid $b$ in $B_{n+1}$ is called {\it quasipositive} if there are braids $w_1,\dots,w_k\in B_{n+1}$ such that $$(\Upsilon ):
b = (w_1 \sigma_{i_1} w_1^{-1})(w_2 \sigma_{i_2} w_2^{-1})\dots(w_k \sigma_{i_k} w_k^{-1})$$
each $w_j$ can be supposed positive.

The purpose of this paper is to solve {\it the quasipositivity problem} in $B_{n+1}$, i.e. to find an algorithm which decides
if a given braid is quasipositive or not.This result can be used in the 16th Hilbert problem, which deals with the classification of real algebraic curves.

We refer the reader to {\it Orevkov [7],[8], Rudolph [9],[10]}
 for more motivations. Our algorithm is based on the ideas of Garside, see {\it Garside [5]}.

 It is based on Theorem 1 below which was conjectured by S.~Orevkov.

Let us denote by $B_{n+1}^+$ the {\it monoid of positive braids}, this is the monoid, (semigroup with the unit) defined by the
 presentation ($\Theta$),see {\it Birman [1], Garside [5]} for details.

It is known by {\it Garside}
that $B_{n+1}^+$ coincides with the submonoid of $B_{n+1}$ generated by $\sigma_1,\dots,\sigma_n$.

This means that two positive words in $\sigma_1,\dots,\sigma_n$ represent the same element of $B_{n+1}$ if and only if one can pass
 from one to another using only {\sl positive} relations.

Moreover theses {\sl positive} relations don't change the number of generators so we can define on $B_{n+1}^+$ a unique length function giving this  number.

 We say that $a$ (left) divides $b$ ($a$ right divides $b$) and we note $a \prec b$ ($b \succ a$)if there exists $c$ satisfying $b=ac$ ($b=ca$).

Let $\Delta$ be {\it the Garside element} of $B_{n+1}$:
$$\Delta = (\sigma_1\sigma_2\dots\sigma_n)(\sigma_1\sigma_2\dots\sigma_{n-1})\dots(\sigma_1\sigma_2)\sigma_1$$

This element is "almost central", i.e.
$\Delta \sigma_k = \sigma_{n+1-k}\Delta$ for $1\leq{{k} \leq{n}}$.Let us denote $\sigma_{n+1-k}$ by $\bar{\sigma}_k$.
We then obtain, for any word $a$ in $B_{n+1}$:  if $s$ even  $\Delta ^{s}a = a\Delta ^{s}$ and if $s$ odd,
$\Delta ^{s}a$ =$ \bar{a}\Delta ^{s}$; in the next,we note $\Delta^s a=\tilde{a} \Delta^s$ when we don't know if $s$ is even or not.

This implies that a positive braid $b$ is right divisible by $\Delta$ in $B_{n+1}^+$ if and only if it is (left) divisible and
 we shall say in this case that $b$ {\it is divisible by} $\Delta$.We call by simple element, every positive word which divides $\Delta$.

Any braid can be represented in the form $\Delta^{-r} w$ with an integer $r$ and a positive braid $w$ which is not divisible by $\Delta$.

Moreover, there is an effective algorithm which allows to compute such a representation, see {\it Garside [5]},{\it Morton, El Rifai [6]}
 and we give here a Garside's lemma
useful for the next.

\newtheorem{lemme0}{Lemma[Garside]0}
\begin{lemme0}

We give positive words $X$,$Y$ such as $\sigma_i X=\sigma_j Y$ for generators $\sigma_i,\sigma_j$.

If $i=j$,then $X=Y$.

If $|i-j|=1$, then there exists $Z$ such as $X=\sigma_j \sigma_i Z$ and $Y=\sigma_i \sigma_j Z$.

If $2 \leq |i-j|$, then there exists $Z$ such as $X=\sigma_j Z$ and $Y=\sigma_i Z$.

\end{lemme0}

In the next paragraph, we give series of results and in the following one their proofs .In the last one we give a computational algorithm and we study its complexity.

\paragraph{Statement of the result}

\newtheorem{definition}{Definition 1}
\begin{definition}
Let $w \in B_{n+1}^+$, $w = u\sigma v$.
We note $$ w =u\sigma v\stackrel{1}{\rightarrow} uv $$ to say that we pass from $w$ to $uv$ by eliminating one generator.This last is said to be a 1-sequence.
So for any $X$ in $B_{n+1}^+$,the k-sequence $$X\stackrel{k}{\rightarrow}Y$$ means that we pass from $X$ to $Y$ by eliminating $k$ generators,i.e there
 exists
$k+1$ positive words $L_j $ with $0\leq{j\leq{k}}$ such as $$L_k = X$$$$L_{j+1}\stackrel{1}{\rightarrow}L_j$$$$L_0=Y$$

\end{definition}

REMARK:

$X\stackrel{k}{\rightarrow}Y$ represents one of the way to reduce $X$ to $Y$.When we have to be more precise, we write

$X=X_k\delta_{i_k} Y_k\stackrel{1}{\rightarrow}X_kY_k=X_{k-1}\delta_{i_{k-1}} Y_{k-1}\stackrel{k-2}{\rightarrow}X_1\delta_{i_1} Y_1
\stackrel{1}{\rightarrow}X_1Y_1=Y$.

Moreover, if there exist others ($X_j', Y_j', \delta_{i_j}'$) reducing $X$ to $Y$, we say that the two k-sequences are equivalent.

\vspace{0.1in}

EXAMPLE:

We note $1$ for $\sigma_1$, $1^{-1}$ for $\sigma_1 ^{-1}$, $2$ for $\sigma_2 $ and so on...

Let  the word $2131$ be in $B_4$ then

$$2132\stackrel{1}{\rightarrow} 212 =121\stackrel{1}{\rightarrow}11$$

that is
$$2132 \stackrel{2}{\rightarrow} 11$$

We see in this example that we can't delete  two generators straight away of the initial word; we have 
to transform $212$ into $121$.

The k-sequence $X\stackrel{k}{\rightarrow}Y$ is  more general than deleting $k$ generators at the same time and in the same
word $X$.

\vspace{0.1in}

We can summarise the results of the whole paper, by:

\newtheorem{lemme}{Fundamental Result}
\begin{lemme}

For positive words  $X$, $Y$, with length$(X)$ = length$(Y) + k$, the following are equivalent:

1. \hspace{0.1in} $Y^{-1}X$ is quasipositive

2. \hspace{0.1in} $X\stackrel{k}{\rightarrow}Y$

3. \hspace{0.1in} $Xd_1d_2d_3...d_k = Yd_1\sigma_{i_1}d_2\sigma_{i_2}d_3...\sigma_{i_{k-1}}d_k\sigma_{i_k}$
for positive $d_j$ and generators $\sigma_{i_j}$

\end{lemme}

The $d_j$'s are defined in the next.
In this case, one quasipositive form of $Y^{-1}X$ is $$Y^{-1}X=
(d_1d_2...d_k)^{-1}d_1\sigma_{i_1}d_2\sigma_{i_2}d_3...\sigma_{i_{k-1}}d_k\sigma_{i_k}$$
Moreover, we have a finite number of quasipositive forms for $Y^{-1}X$; each form corresponds to a way to reduce $X$ to $Y$, by
$X\stackrel{k}{\rightarrow}Y$, and theses last are finite.

\newtheorem{lemme1}{Lemma 1}
\begin{lemme1}

For positive words $A$, $X$, $Y$, and a positive integer $k$,the k-sequence
$AX \stackrel{k}{\rightarrow} AY$ is equivalent to  the other one  $X \stackrel{k}{\rightarrow} Y$

\end{lemme1}

\newtheorem{guess1}{Theorem1}
\begin{guess1}

We give $b=\Delta^{-r}b_1...b_s = N^{-1}P$ written in its normal form, with algebraic length $k$, then
$b$ is quasipositive if and only if $P \stackrel{k}{\rightarrow} N$.

\end{guess1}

PROOF:

Let  $b$ be quasipositive,  $b=\Delta^{-r}b_+$ :  by definition $b=w_1\sigma_{i_1}{w_1}^{-1}...w_k\sigma_{i_k}{w_k}^{-1}$ for positive $w_i$ . 
Each $w_i$ written in its normal form gives $w_i = \Delta^{m_i}{w_i}_+$  and also
${w_i}^{-1}=\Delta^{-p_i}{v_i}_+$ for $m_i,p_i$ positive integer.
Moreover $w_i{w_i}^{-1} =1$ is equivalent to  $\Delta^{m_i}{w_i}_+\Delta^{-p_i}{v_i}_+ =1$, $\Delta^{m_i-p_i}\tilde{{w_i}_+}{v_i}_+=1$; finally 
$\tilde{{w_i}_+}{v_i}_+=\Delta^{p_i-m_i}$.

\vspace{0.2in}

We can write $b$ as  $b=\Delta^{m_1-p_1+...+m_k-p_k} \tilde{{w_1}_+}\tilde{\sigma_{i_1}} {v_1}_+...\tilde{{w_k}_+}\tilde{\sigma_{i_k}} {v_k}_+$, and we note $w$ for $\tilde{{w_1}_+}\tilde{\sigma_{i_1}} {v_1}_+...\tilde{{w_k}_+}\tilde{\sigma_{i_k}} {v_k}_+$ and $m-p$ for ${m_1-p_1+...+m_k-p_k}$.

If we cancel in $w$ all the generators $\sigma_{i_j}$,
$w$ becomes $\tilde{{w_1}_+}{v_1}_+...\tilde{{w_k}_+}{v_k}_+$ which is equal to $\Delta^{p_1-m_1+...+p_k-m_k}=\Delta^{p-m}$.In other words,

$$w \stackrel{k}{\rightarrow}\Delta^{p_1-m_1+...+p_k-m_k} $$

Futhermore, $b$ in its normal form satisfies $w=\Delta^u b_+$, for an integer $u$.

As $b=\Delta^{-r} b_+$ and $b=\Delta^{m-p}w$, we obtain$\Delta^{m-p}w =\Delta^{-r} b_+$, $w=\Delta^{p-m-r}b_+$.

Finally $r+u=p-m$ and the sequence becomes

$$w=\Delta^u b_+ \stackrel{k}{\rightarrow} \Delta^{p-m}=\Delta^u \Delta^r$$

By  lemmas above, $AX \stackrel{k}{\rightarrow} AY$ beeing equivalent to  $X \stackrel{k}{\rightarrow} Y$,
 we obtain here, by simplifying by $\Delta^u$,
$$b_+ \stackrel{k}{\rightarrow}\Delta^r$$

Moreover $b=N^{-1}P$ satisfies by definition $RP=b_+$ and $RN =\Delta^r$ for a positive $R$: the sequence above becomes

$b_+=RP \stackrel{k}{\rightarrow}\Delta^r=RN$ which is equivalent to $$P \stackrel{k}{\rightarrow} N$$

\vspace{0.2in}

Conversely, if  $b=N^{-1}P$ satisfies   $P \stackrel{k}{\rightarrow} N$, we can write $P={\sigma_{i_1}}...{\sigma_{i_k}}N$ and obviously
 $b=N^{-1}{\sigma_{i_1}}...{\sigma_{i_k}}N$ is quasipositive.

\vspace{0.2in}

EXAMPLE:

In $B_4$, $b =212^{-1}2^{-1}32$ is quasipositive, with $k=2$: $b =1^{-1}1^{-1}2132$ in its normal form; we have

$P= 2132 \stackrel{2}{\rightarrow} N=11$; indeed

$$2132 \stackrel{1}{\rightarrow} 212=121\stackrel{1}{\rightarrow}11$$

\vspace{0.3in}

Theorem 11 gives us a therotical way to decide whether a word is quasipositive or not.We now give some results, proven in the following,
 which lead us to
a computational algorithm.

Dealing with a k-sequence $X \stackrel{k}{\rightarrow} Y$, we saw that the $k$ generators deleted are not necessary deleted at the same time;
but we see in the next that we can trasform any sequence into one of the form

$$A_1\sigma_{i_1}A_2...A_k\sigma_{i_k}A_{k+1} \stackrel{k}{\rightarrow} A_1A_2...A_kA_{k+1}$$

where the $k$ generators can be deleted in the same word and at the same time.

\newtheorem{lemme2}{Lemma 2}
\begin{lemme2}

$b=\Delta^{-r}b_+$ with algebraic length $k$ is quasipositive if and only if there exist $k+1$ positive words $A_1,...,A_{k+1}$ and $k$ 
generators
$\sigma_{i_1},...,\sigma_{i_k}$ such as

$$b_+=A_1\sigma_{i_1}A_2\sigma_{i_2}...A_k\sigma_{i_k}A_{k+1}$$
$$\Delta^{-r}=A_1A_2...A_kA_{k+1}$$

\end{lemme2}

This lemma gives a stronger result than $b_+\stackrel{k}{\rightarrow}\Delta^{-r} $.Unfortunately, we can't obtain a similar version 
for $P \stackrel{k}{\rightarrow} N$.But we have a slighty different one, which will give the main tool for an algorithm.

\newtheorem{lemme3}{Lemma 3}
\begin{lemme3}

We give $b=\Delta^{-r}b_+=\Delta^{-r}b_1b_2...b_r...b_s$ with algebraic length $k$. 
We define for $0 \leq j\leq r$,

$$P_j=b_{j+1}...b_{s-r+j}$$
 and $N_j$ such as 

$$b_1...b_jN_jb_{s-r+j+1}...b_s=\Delta^{r}$$

$b$ is quasipositive if and only if for one $j$, $0 \leq j\leq r$, there exist $k$ positive words $A_2,...,A_{k+1}$ (here $A_1 = 1$)
and $k$ generators $\sigma_{i_1},...,\sigma_{i_k}$ such as 

$$P_j=\sigma_{i_1}A_2\sigma_{i_2}...A_k\sigma_{i_k}A_{k+1}$$
$$N_j=A_2...A_kA_{k+1}$$

\end{lemme3}

In this case, we output the quasipositive form of $b$:we set $T_j=(b_{s-r+j+1}...b_s)$ and

$T_j b T_j^{-1} = N_j^{-1}P_j$ so $$b = T_j^{-1}[(A_2...A_kA_{k+1})^{-1}(\sigma_{i_1}A_2\sigma_{i_2}...A_k\sigma_{i_k}A_{k+1})]T_j$$


\paragraph{Technical proofs}

The next lemma is a classical known result.

\newtheorem{lemme4}{Lemma 4}
\begin{lemme4}
lgcd (left greatest common divisor)is distributive.

Let $R$,$A$,$B$ be words in $B_{n+1}$,lgcd $(RA,RB)=R$ lgcd $(A,B)$.

(resp. rgcd (right greatest common divisor) is distributive, rgcd $(AR,BR)=$ rgcd $(A,B)R$).

\end{lemme4}

PROOF:

$A=$ lgcd $(A,B)A'$ and $B=$ lgcd $(A,B)B'$ by definition of lgcd, for positive $A'$,$B'$.

So $RA=R$ lcd $(A,B)A'$ and $RB=R$ lgcd $(A,B)B'$, that is 

$R$ lgcd $(A,B)\prec$ lgcd $(RA,RB)$, lgcd $(RA,RB)$ beeing the left greatest common divisor.

Let $U$ be such as lgcd $(RA,RB)=R$ lgcd $(A,B)U$.

Then lgcd $(RA,RB)=R$ lgcd $(A,B)U\prec RA =R$ lgcd $(A,B)A'$ gives 

$R$ lgcd $(A,B)U \prec R$ lgcd $(A,B)A'$, so
$U\prec A'$ and in the same way $U\prec B'$.

Finally, lgcd $(A,B)U$ divides 
lgcd $(A,B)A'= A$ and also $B$; 

but lgcd $(A,B)$ beeing
  the left greatest common divisor of $A$ and $B$ , $U=1$ and we have the result,

  lgcd $(RA,RB)=R$ lgcd $(A,B)$.

\vspace{0.3in}

The next one, very important for the next, already appeared in some work of Juan Gonzalez-Meneses.

\newtheorem{lemme5}{Lemma 5}
\begin{lemme5}
Let $w$ be a word in $B_{n+1}^+$ and $\sigma$ a generator.

We note ${D_w}^l(\sigma )=$ (X positive such as $wX\prec \sigma wX$).

In this case $wX\sigma'= \sigma wX$ for a generator $\sigma'$.

${D_w}^l(\sigma )$ is not empty and has a minimal element for left divisibility, called  ${d_w}^l(\sigma )$,this last divides
all $X \in {D_w}^l(\sigma)$.

(We have  same results for right side: ${d_w}^r(\sigma )$ minimal for right divisibility in
 ${D_w}^r(\sigma )$ =(Y positive such that $Y w\sigma \succ Yw$)).

\end{lemme5}

PROOF:

Let $w$ be written in its left normal form, $$w=\Delta^{-r}b_1...b_s$$ with $-r$ integer and let $a_i$ be such as $b_ia_i=\Delta$ for
$1\leq i\leq s$.

We have $b_1...b_{s-1}b_s.a_s{\bar{a_{s-1}}}...{\tilde{a_1}}=\Delta^{s}$ and if we note $X=a_s{\bar{a_{s-1}}}...{\tilde{a_1}}$
then $X$ satisfies
$\sigma wX=\sigma \Delta^{-r}b_1...b_s.a_s{\bar{a_{s-1}}}...{\tilde{ba_1}}=\sigma \Delta^{-r+s}=\Delta^{-r+s} {\tilde{\sigma}}$
that is $$\sigma wX=wX{\tilde{\sigma}}$$ and $X \in {D_w}^l(\sigma )$.

Let's now prove that ${D_w}^l(\sigma )$ has a minimal element:
if $X$ and $Z$ are in ${D_w}^l(\sigma )$, also is lgcd$(X,Z)$,indeed
$$\sigma wX=wX\sigma'$$
$$\sigma wZ=wZ\sigma''$$ for generators $\sigma'$, $\sigma''$.

Then lgcd$(\sigma wX,\sigma wZ)=$ lgcd$(wX\sigma',wZ\sigma'')$ and by Lemma 31 on the distributivity of lgcd,
$\sigma w$ lgcd$(X,Z)=w$ lgcd$(X\sigma',Z\sigma'')$.

But  lgcd$(X,Z)\prec X \prec X \sigma'$ and lgcd$(X,Z)\prec Z \prec Z \sigma''$ so

 lgcd $(X,Z) \prec $ lgcd $(X\sigma' ,Z\sigma'')$ and $w$ lgcd $(X,Y)u = w$ lgcd $(X\sigma',Z\sigma '')$ for a positive $u$.

If we replace this in what before,$\sigma w$lgcd$(X,Z)=w$lgcd$(X,Z)u$.

Now, dealing with positive words, the algebraic length of the two factors in the equality must be equal, so $u$ is a generator, $u=\beta$

Finally $\sigma w$ lgcd $(X,Z)=w$ lgcd $(X,Z)u$ gives $\sigma w$ lgcd $(X,Z)=w$ lgcd $(X,Z)\beta$ so lgcd $(X,Z)$ satisfies the 
same properties as $X,Z$,

lgcd $(X,Y) \in {D_w}^l(\sigma )$.

And ${d_w}^l(\sigma )$ can be  taken as beeing the commun lgcd of all  words in ${D_w}^l(\sigma )$.

\newtheorem{lemme 6}{Lemma 6}
\begin{lemme 6}
Let $w$ be a positive word, $\sigma$ a generator and ${d_w}^l(\sigma )$ as above:

$\sigma w{d_w}^l(\sigma )= w{d_w}^l(\sigma )\sigma'$.

If $f$ and $h$ are positive words such as $$\sigma wf=wh$$ then, ${d_w}^l(\sigma )\prec f$,
${d_w}^l(\sigma )\prec h$ and
if $f={d_w}^l(\sigma )x$,
 $h={d_w}^l(\sigma )\sigma' x$ for a positive $x$.

(We have same results for right divisibility).

\end{lemme 6}

PROOF:

We use Garside's Lemma to prove it when $w$ is a generator: $w=\beta$

if $\sigma=\beta$ then $\sigma\beta=\beta\sigma$ and ${d_w}^l(\sigma ) = 1$, ${d_w}^l(\sigma )$ obviously divides $f$  as above

if $\sigma$ and $\beta$ commutes then again $\sigma\beta=\beta\sigma$ and ${d_w}^l(\sigma ) = 1$ as above

if $\sigma\beta\sigma=\beta\sigma\beta$, then ${d_w}^l(\sigma ) = \sigma$ and we know:

 $\sigma\beta f=\beta h$ gives $f=\sigma x$ and $h=\sigma \beta x$ for a positive $x$ which gives 

 ${d_w}^l(\sigma ) \prec (f, h)$, $f={d_w}^l(\sigma )x$, $h={d_w}^l(\sigma ) \sigma' x$, $ \sigma'= \beta$.

(we have same results for right divisibility)

So let's write now the proof in the general case:$$\sigma w{d_w}^l(\sigma )=w{d_w}^l(\sigma ) \sigma'$$

                                         $$\sigma wf=wh$$

We have lgcd$(\sigma w{d_w}^l(\sigma ),\sigma wf)=$lgcd$(w{d_w}^l(\sigma ) \sigma',wh)$ which is equivalent to
$\sigma w$lgcd$({d_w}^l(\sigma ),f)=w$lgcd$({d_w}^l(\sigma )\sigma',h)$.

If we note $a=$lgcd$({d_w}^l(\sigma ),f)$ and $b=$lgcd$({d_w}^l(\sigma)\sigma',h)$:

$au={d_w}^l(\sigma )$ for a positive $u$ and $\sigma wa=wb$ gives $\sigma wau=wbu$, that is $\sigma w{d_w}^l(\sigma )=wbu$ 
but $\sigma w{d_w}^l(\sigma )=w{d_w}^l(\sigma )\sigma'$ so $$bu={d_w}^l(\sigma )\sigma'$$

If $u$ no equals to $1$ then $u=v\delta$ for a generator $\delta$:

 $au={d_w}^l(\sigma)$ gives $au\sigma'=bu$ that is $av\delta\sigma'=bv\delta$.

The word $\delta$  'pseudocommutes' with $\sigma'$ and  ${d_{\delta}}^r(\sigma')$ satisfies ${d_{\delta}}^r(\sigma') \delta \sigma'
= \sigma'' {d_{\delta}}^r(\sigma')\delta$; there exists $y$ positive, as a consequence of Lemma 5 1, with

$$av=y{d_{\delta}}^r(\sigma')$$
$$bv=y \sigma'' {d_{\delta}}^r(\sigma')$$

Finally,$$\sigma w{d_w}^l(\sigma) = w{d_w}^l(\sigma ) \sigma'$$

$$\sigma wau=wbu$$

$$\sigma wav\delta=wbv\delta$$

$$\sigma wy{d_{\delta}}^r(\sigma')\delta=wy\sigma''{d_{\delta}}^r(\sigma' )\delta$$

$$\sigma wy=wy\sigma''$$

and $y$ satisfies the same  propriety as ${d_w}^l(\sigma )$, $y$ belongs to ${D_w}^l(\sigma )$.

But $y{d_{\delta}}^r(\sigma')\delta=d$ and $y$ is strictly smaller than  ${d_w}^l(\sigma)$;
which is a contradiction  of the minimality of ${d_w}^l$ in ${D_w}^l(\sigma )$.

As all this results from $u$  not equal to $1$ ,we obtain $u=1$ and

lgcd$({d_w}^l(\sigma),f)=a=au={d_w}^l(\sigma)$

lgcd$({d_w}^l(\sigma )\sigma',h)=b=bu={d_w}^l(\sigma )\sigma'$

so ${d_w}^l(\sigma )\prec f$,${d_w}^l(\sigma )\sigma'\prec h$;

if we put $f={d_w}^l(\sigma )x$,then $\sigma wf=wh$ gives $$\sigma w{d_w}^l(\sigma )x=wh$$
$$w{d_w}^l(\sigma )\sigma'x=wh$$ and $h={d_w}^l(\sigma )\sigma'x$, which is what we are looking for.

\newtheorem{lemme7}{Lemma 7}
\begin{lemme7}

Let  $X \stackrel{k}{\rightarrow} Y$ be a k-sequence, then there exits $k$ positive words $d_1,...,d_k$ and $k$ 
generators $\sigma_{i_1},...,\sigma_{i_k}$
such as
$$Xd_1d_2...d_k =Yd_1\sigma_{i_1}d_2\sigma_{i_2}...d_k\sigma_{i_k}$$

\end{lemme7}

PROOF:

 $X \stackrel{k}{\rightarrow} Y$ means there exist $X_j,Y_j,\delta_{i_j}$ ,$1\leq j \leq k$
satisfying 

$$X=X_k\delta_{i_k}Y_k\stackrel{1}{\rightarrow}X_kY_k=X_{k-1}\delta_{i_{k-1}}Y_{k-1}\stackrel{k-2}{\rightarrow}
X_2Y_2=X_1\delta_{i_1}Y_1\stackrel{1}{\rightarrow}Y=X_1Y_1$$

We now consider :

$d_1={d_{Y_1}}^r(\delta_{i_1})$ which satisfies by Lemma 5 1, $\delta_{i_1}Y_1d_1=Y_1d_1\delta_{i_1}'$

$d_2={d_{Y_2d_1}}^r(\delta_{i_2})$ such as $\delta_{i_2}{Y_2d_1}d_2={Y_2d_1}d_2\delta_{i_2}'$

$d_j={d_{Y_jd_1d_2...d_{j-1}}}^r(\delta_{i_j})$ such as $\delta_{i_j}{Y_jd_1d_2...d_{j-1}}d_j={Y_jd_1d_2...d_{j-1}}d_j\delta_{i_j}'$

$d_k={d_{Y_kd_1d_2...d_{k-1}}}^r(\delta_{i_k})$ such as $\delta_{i_k}{Y_kd_1d_2...d_{k-1}}d_k={Y_jkd_1d_2...d_{k-1}}d_k\delta_{i_k}'$

By definition of theses $d_j$, we have by multiplying the sequence above
 by $d_1d_2...d_k$:

for the first line , $Yd_1d_2...d_k$

for the second one, $X_2Y_2d_1d_2...d_k=X_1\delta_{i_1}Y_1d_1d_2...d_k=X_1Y_1d_1\delta_{i_1}'d_2...d_k=Yd_1\delta_{i_1}'d_2...d_k$; 
moreover this equality gives $X_2Y_2d_1=Yd_1\delta_{i_1}'$

for the third one, $X_3Y_3d_1d_2...d_k=X_2\delta_{i_2}Y_2d_1d_2...d_k=X_2Y_2d_1d_2\delta_{i_2}'d_3...d_k=Yd_1\delta_{i_1}'d_2\delta_{i_2}'
d_3...d_k$ by what above; and we also obtain $X_3Y_3d_1d_2=Yd_1\delta_{i_1}'d_2\delta_{i_2}'$

for the $k$th one, $X_kY_kd_1d_2...d_k=X_{k-1}\delta_{i_{k-1}}Y_{k-1}d_1d_2...d_k=X_{k-1}Y_{k-1}d_1d_2...d_{k-1}\delta_{i_{k-1}}'d_k=Yd_1
\delta_{i_1}'d_2\delta_{i_2}'d_3\delta_{i_3}'...d_{k-1}\delta_{i_{k-1}}'d_k$
by induction,

and in the same way, the final result $Xd_1d_2...d_k=Yd_1\delta_{i_1}'d_2\delta_{i_2}'d_3\delta_{i_3}'...d_{k-1}\delta_{i_{k-1}}'d_k
\delta_{i_k}'$

\vspace{0.2in}

In the next Lemma, we see that within the condition 

$Xd_1d_2...d_k =Yd_1\sigma_{i_1}d_2\sigma_{i_2}...d_k\sigma_{i_k}$, our $d_j$ are minimal.

\newtheorem{lemme8}{Lemma 8}
\begin{lemme8}

Let  $X=X_k\delta_{i_k} Y_k \stackrel{1}{\rightarrow} X_kY_k=X_{k-1}\delta_{i_{k-1}} Y_{k-1}\stackrel{k-1}{\rightarrow}X_1Y_1=Y$
be a k-sequence.

 If there exits $k$ positive words $A_1,...,A_k$ and $k$ generators $\sigma_{i_1},...,\sigma_{i_k}$
such as
$XA_1A_2...A_k =YA_1\sigma_{i_1}A_2\sigma_{i_2}...A_k\sigma_{i_k}$, then we can find an equivalent k-sequence

$X=X_k'\delta_{i_k}' Y_k'=L_k \stackrel{1}{\rightarrow} X_k'Y_k'=X_{k-1}'\delta_{i_{k-1}}' Y_{k-1}'=L_{k-1}
\stackrel{k-1}{\rightarrow}X_1'Y_1'=Y=L_0$ such as

$$L_jA_1...A_j=YA_1\sigma_{i_1}A_2\sigma_{i_2}...A_j\sigma_{i_j}$$

\end{lemme8}

PROOF:

By induction on $l$ = length $(Y)$.

\vspace{0.1in}

$l=0$.Which is equivalent to $Y=1$. Let   $X=\stackrel{k}{\rightarrow} 1$ be a k-sequence with
$XA_1A_2...A_k =A_1\sigma_{i_1}A_2\sigma_{i_2}...A_k\sigma_{i_k}$.We use now an induction on $k$.

\hspace{0.1in} $k=1$ gives $X=X_1\delta_{i_1} Y_1\stackrel{1}{\rightarrow}1$ and $XA_1=A_1\sigma_{i_1}$.
Obviously $L_1=X$ satisfies $L_1A_1=A_1\sigma_{i_1}$, so the result.

\vspace{0.1in}

\hspace{0.1in} We suppose the resut true for any (k-1)-sequence.With the previous datas,we set 
$$U=A_1\sigma_{i_1}A_2\sigma_{i_2}...A_{k-1}\sigma_{i_{k-1}}.(A_1...A_{k-1})^(-1)$$

Which gives $UA_1...A_{k-1}=YA_1\sigma_{i_1}A_2\sigma_{i_2}...A_{k-1}\sigma_{i_{k-1}}$ and by replacing in
$XA_1A_2...A_k =A_1\sigma_{i_1}A_2\sigma_{i_2}...A_k\sigma_{i_k}$, we obtain $XA_1A_2...A_k=UA_1A_2...A_k\sigma_{i_k}$.

$U$ is positive and by using $D$ minimal, as in Lemma 5 1, such as $DA_1A_2...A_k\sigma_{i_k}= \sigma_{i_k}'DA_1A_2...A_k$,

$$X=E \sigma_{i_k}' D$$
$$U=ED$$

for $E$ positive.

Moreover,length $(U)=k-1$ gives $U \stackrel{k-1}{\rightarrow}1$. As,
$UA_1...A_{k-1}=YA_1\sigma_{i_1}A_2\sigma_{i_2}...A_{k-1}\sigma_{i_{k-1}}$, the hypothesis of induction gives

$$U=L_{k-1}\stackrel{k-j}{\rightarrow}L_j\stackrel{j}{\rightarrow}L_0=1$$

with $L_jA_1...A_j=A_1\sigma_{i_1}A_2\sigma_{i_2}...A_j\sigma_{i_j}$.If we add
$X=E \sigma_{i_k}' D\stackrel{1}{\rightarrow}U=ED$, we obtain a k-sequence with each line $L_j$ satisfying the final result.

\vspace{0.1in}

For any $l$.We suppose the result true for $l-1$.Let   $X=\stackrel{k}{\rightarrow} Y$ with length $(Y)=l$.

With $Y=Y'\delta$, we obtain a (k+1)-sequence

$$X=\stackrel{k}{\rightarrow} Y=Y'\delta \stackrel{1}{\rightarrow}Y'$$

By setting $1=A_0$ and $\delta =\sigma_{i_0}$, the equality $XA_1A_2...A_k =YA_1\sigma_{i_1}A_2\sigma_{i_2}...A_k\sigma_{i_k}$ becomes
$XA_0A_1A_2...A_k =Y'A_0\sigma_{i_0}A_1\sigma_{i_1}A_2\sigma_{i_2}...A_k\sigma_{i_k}$.

As length $(Y')=l-1$, we have by hypothesis of induction, an equivalent (k+1)-sequence
$X=L_k\stackrel{k-j}{\rightarrow}L_j \stackrel{j}{\rightarrow}L_0 \stackrel{1}{\rightarrow}L_0'=Y'$.
Each $L_j$ satisfies $L_jA_0A_1...A_j=Y'A_0\sigma_{i_0}A_1\sigma_{i_1}A_2\sigma_{i_2}...A_j\sigma_{i_j}$.
And $L_0$ is such as $L_0A_0=Y'A_0\sigma_{i_0}$
which is $L_0=Y$.

And the final result $X=L_k\stackrel{k-j}{\rightarrow}L_j \stackrel{j}{\rightarrow}L_0 =Y$, with
$L_jA_0A_1...A_j=Y'A_0\sigma_{i_0}A_1\sigma_{i_1}A_2\sigma_{i_2}...A_j\sigma_{i_j}$ eqivalent to
$L_jA_1...A_j=YA_1\sigma_{i_1}A_2\sigma_{i_2}...A_j\sigma_{i_j}$.

\newtheorem{lemme9}{Lemma 9}
\begin{lemme9}

Conversely, if there exits $k$ positive words $A_1,...,A_k$ and $k$ generators $\sigma_{i_1},...,\sigma_{i_k}$
such as
$XA_1A_2...A_k =YA_1\sigma_{i_1}A_2\sigma_{i_2}...A_k\sigma_{i_k}$, then

 $$X \stackrel{k}{\rightarrow} Y$$

\end{lemme9}

PROOF:

By induction on $l=$length$(Y)$.

\vspace{0.1in}

$l=0$: which is equivalent to $Y=1$.

$XA_1A_2...A_k =YA_1\sigma_{i_1}A_2\sigma_{i_2}...A_k\sigma_{i_k}$ becomes 
$XA_1A_2...A_k =A_1\sigma_{i_1}A_2\sigma_{i_2}...A_k\sigma_{i_k}$.The equality of length in the positive monoid gives 
$l(X)=k$ and obviously
$$X \stackrel{k}{\rightarrow} 1$$

\vspace{0.2in}

Any $l$. For $Y'$ with length $l-1$ we suppose the result  true: for any $k$ integer and
$XA_1A_2...A_k =YA_1\sigma_{i_1}A_2\sigma_{i_2}...A_k\sigma_{i_k}$ then  $X \stackrel{k}{\rightarrow} Y$.
We give now $Y=Y'\delta $ of length $l$ satisfying
$XA_1A_2...A_k =Y' \delta A_1\sigma_{i_1}A_2\sigma_{i_2}...A_k\sigma_{i_k}$ for positive $A_i$'s.

We consider $A_0=1$ and $\sigma_{i_0}=\delta$ to have a new equality 
$XA_0A_1A_2...A_k =Y' A_0\sigma_{i_0} A_1\sigma_{i_1}A_2\sigma_{i_2}...A_k\sigma_{i_k}$.

But $l(Y')=l-1$, we can apply what before to obtain $X \stackrel{k+1}{\rightarrow} Y'$.
Lemma 8 1 tells us that we have an equivalent k-sequence such as

$X=L_k\stackrel{k-j}{\rightarrow}L_j\stackrel{j}{\rightarrow}L_0\stackrel{1}{\rightarrow}L_0'=Y'$
and $L_jA_0A_1...A_j=Y'A_0\sigma_{i_0}A_1\sigma_{i_1}A_2\sigma_{i_2}...A_j\sigma_{i_j}$.
In a particular case, $L_0A_0=Y'A_0\sigma_{i_0}$ which is $L_0=Y$.

And the final k-sequence $X \stackrel{k}{\rightarrow} L_0=Y$.

We now give the proof of Lemma 1 1.

\vspace{0.1in}

PROOF Lemma 1 1:

If $X\stackrel{k}{\rightarrow}Y$, obviously for any positive $A$, $AX\stackrel{k}{\rightarrow}AY$.

\vspace{0.1in}

Conversely if we have $AX\stackrel{k}{\rightarrow}AY$, Lemma 7 1 gives $d_1,d_2,...,d_k$ and $\sigma_{i_1},\sigma_{i_2},...,\sigma_{i_k}$ such as

$$AXd_1d_2...d_k =AYd_1\sigma_{i_1}d_2\sigma_{i_2}...d_k\sigma_{i_k}$$

This last equality is equivalent to $Xd_1d_2...d_k =Yd_1\sigma_{i_1}d_2\sigma_{i_2}...d_k\sigma_{i_k}$ when simplifying by $A$.

And Lemma 9 1 finally gives $X\stackrel{k}{\rightarrow}Y$.

\vspace{0.2in}

The next two lemmas will lead to the proof of Lemma 2 1.

\newtheorem{lemme10}{Lemma 10}
\begin{lemme10}

Let  $A\sigma B \stackrel{1}{\rightarrow}AB$ be a 1-sequence for positive $A,B$.

We note the  line $A\sigma B$ by $L_1$ and the other one $AB$ by $L_0$.

Then there exist positive $A',B'$, a generator $\sigma'$ such as 

$$A\sigma B =A'\sigma' B'$$
$$AB=A'B'$$
$$L_1=A\sigma B =A'\sigma' B'\stackrel{1}{\rightarrow}L_0=AB=A'B'$$
 and rgcd $(L_1,L_0) =B'$.

(In the same way, we have $A\sigma B =A''\sigma'' B''$ and 
$AB=A''B''$ for positive $A'',B''$ such as $A''=$lgcd $(L_1,L_0)$).

\end{lemme10}

PROOF:

We note rgcd$(L_1,L_0)=r$ and obviously $r\succ B$: $r=uB$ for a positive $u$.

We also have $A\sigma=xu$, $A=yu$ for positive $x,y$ and 

$$yu\sigma =xu$$

By Lemma 3 1, $y=y'd$ and $x=x'\sigma'd$ with $d={d_u}^r(\sigma)$ :

$L_1=A\sigma B =x'\sigma'duB$ and $L_0=AB=y'duB$.

Theses last give  $(L_1,L_0) \succ duB$ and $uB=r \succ duB$.

Finally $d=1$ and we put $A'=y'$,$B'=uB$ to obtain the result.(It's the same for lgcd).

\vspace{0.2in}

In a more general case, we obtain the next result.

\newtheorem{lemme11}{Lemma 11}
\begin{lemme11}

 Let $L_k = A_1\sigma_{i_1}A_2\sigma_{i_2}A_3...A_k\sigma_{i_k}A_{k+1}\stackrel{1}{\rightarrow}
L_{k-1} = A_1\sigma_{i_1}A_2\sigma_{i_2}A_3...A_kA_{k+1}\stackrel{k-2}{\rightarrow}L_1 = A_1\sigma_{i_1}A_2...A_kA_{k+1}\stackrel{1}{\rightarrow} L_0 = A_1A_2...A_kA_{k+1}$ be a k-sequence.

Each line $L_j$, $0 \leq j \leq k$, satisfies $L_j= A_1\sigma_{i_1}A_2\sigma_{i_2}...\sigma_{i_{j-1}}A_j\sigma_{i_j}A_{j+1}...A_kA_{k+1}$ where $\sigma_{i_r}$ for $j+1 \leq r \leq k$ are cancelled.

This sequence can be written in an equivalent way,

 $L_k = A_1'\sigma_{i_1}'A_2'\sigma_{i_2}'A_3'...A_k'\sigma_{i_k}'A_{k+1}'\stackrel{1}{\rightarrow}L_{k-1} = A_1'\sigma_{i_1}'A_2'
 \sigma_{i_2}'A_3'...A_k'A_{k+1}'
 \stackrel{k-2}{\rightarrow}L_1 = A_1'\sigma_{i_1}'A_2'...A_k'A_{k+1}'\stackrel{1}{\rightarrow} L_0 = A_1'A_2'...A_k'A_{k+1}'$

  such as rgcd $(L_k,...,L_1,L_0)=A'_{k+1}$.

\end{lemme11}

PROOF:

 Let $L_k = A_1\sigma_{i_1}A_2\sigma_{i_2}A_3...A_k\sigma_{i_k}A_{k+1}\stackrel{k}{\rightarrow} L_0 = A_1A_2...A_kA_{k+1}$ be a k-sequence.

We say this last satisfies (E) if $L_k = A_1'\sigma_{i_1}'A_2'\sigma_{i_2}'A_3'...A_k'\sigma_{i_k}'A_{k+1}'\stackrel{k}{\rightarrow} L_0 = A_1'A_2'...A_k'A_{k+1}'$
with rgcd $(L_k,...,L_0)= A_{k+1}'$.

\vspace{0.1in}

We use here a double induction.Iduction on $k$ and on $l=$length $(L_0)$.

E($k$,$l$) is true, means that a k-sequence with length $(L_0)=l$, satisfies (E) as above.

The basis of induction are as follow:

1- E($k'$,$0$) is true for any $k'$. Indeed, any k'-sequence with length $(L_0)=0$ gives $A_1=...=A_{k+1}=1$.Obviously,
rgcd $(L_k,...,L_1,L_0=1)=1=A_{k+1}$.

Here, the $A_i'$'s can be taken equals to $A_i=1$'s.

2- E($1$,$l$) is true for any $l$, by Lemma 10 1.

\vspace{0.1in}

We can then suppose E($k'$,$l-1$) is true for any $k'$ and E($k'$,$l$) is true for any $0 \leq k' \leq k-1$, and let us prove that 

E($k$,$l$) is true.

\vspace{0.1in}

We have the k-sequence $L_k = A_1\sigma_{i_1}A_2\sigma_{i_2}A_3...A_k\sigma_{i_k}A_{k+1}\stackrel{k}{\rightarrow}
 L_0 = A_1A_2...A_kA_{k+1}$.
We note $R=$ rgcd $(L_k,...,L_0)$.

We consider the (k-1)-sequence $L_{k-1} = A_1\sigma_{i_1}A_2\sigma_{i_2}A_3...\sigma_{i_{k-1}}A_kA_{k+1}\stackrel{k-1}{\rightarrow}
 L_0 = A_1A_2...A_kA_{k+1}$.

As E($k-1$,$l$) is true by hypothesis of induction, we have

$$L_{k-1} = A_1'\sigma_{i_1}'A_2'\sigma_{i_2}'A_3'...\sigma_{i_{k-1}}'A_k'\stackrel{k-1}{\rightarrow} L_0 = A_1'A_2'...A_k'$$

with rgcd $(L_{k-1},...,L_0)=A_k'$.

Moreover, $L_j \succ A_kA_{k+1}$ for $0 \leq j \leq k-1$, so rgcd $(L_{k-1},...,L_0)=A_k' \succ A_kA_{k+1}$. And $A_k'=VA_kA_{k+1}$ for

a positive $V$.

We can  finally write

$L_k=A_1'\sigma_{i_1}'A_2'\sigma_{i_2}'A_3'...\sigma_{i_{k-1}}'VA_k'\sigma_{i_k}A_{k+1}$.

Now, if we write $A_1'=\sigma B$, we obtain the following {k+1}-sequence

$L_k=\sigma B\sigma_{i_1}'A_2'\sigma_{i_2}'A_3'...\sigma_{i_{k-1}}'VA_k\sigma_{i_k}A_{k+1}\stackrel{k}{\rightarrow}
L_0=\sigma BA_2'...A_{k-1}'VA_kA_{k+1}\stackrel{1}{\rightarrow}L_0'= BA_2'...A_{k-1}'VA_kA_{k+1}$.

$R=$ rgcd $(L_k,...,L_0)$ is such as $A_k'=VA_kA_{k+1}=$ rgcd $(L_k,...,L_0) \succ R$, so $L_0' \succ R$ and $R$ 
common to $L_k,...,L_0,L_0'$

satisfies rgcd $(L_k,...,L_0,L_0')=R$.

\vspace{0.1in}

We deal wih a {k+1}-sequence such as length$(L_0')=l-1$; as E($k+1$,$l-1$) is true, this lat sequence can be written,

$L_k=A_0''\sigma'' A_1''\sigma_{i_1}''A_2''\sigma_{i_2}''A_3''...\sigma_{i_{k-1}}''A_k''\sigma_{i_k}''R
\stackrel{k}{\rightarrow}L_0=A_0''\sigma'' A_1''A_2''...A_{k-1}''A_k''R
\stackrel{1}{\rightarrow}L_0'= A_0''A_1''A_2'...A_{k-1}''A_k''R$.

Finally, by considering $L_k\stackrel{k}{\rightarrow}L_0$ with the above $A_i''$'s, our initial sequence is written in the right way.

\vspace{0.2in}

We give now the proof of Lemma 2 1.

\vspace{0.1in}

PROOF Lemma 2 1:

For $b=\Delta^{-r}b_+$ quasipositive, with algebraic length $k$, we have the k-sequence $b_+\stackrel{k}{\rightarrow}\Delta^r$.

Lemma 7 1 gives $d_1,d_2,...,d_k$ such as

$b_+d_1d_2...d_k=\Delta^r d_1\sigma_{i_1}d_2\sigma_{i_2}...d_k\sigma_{i_k}$.

Let $s$ be an integer as large as needed such as $d_1d_2...d_k \prec \Delta^{s}$ and $r \leq s$.

We have, with $d_1d_2...d_kU=\Delta^s$, for a positive $U$, 
$b_+\Delta^s =b_+d_1d_2...d_kU=\Delta^r d_1\sigma_{i_1}d_2\sigma_{i_2}...d_k\sigma_{i_k}U$
by Lemma 7 1.

We can suppose $s$ even, and $b_+\Delta^s =\Delta^s b_+ =\Delta^r \Delta^{s-r}b_+$ and,

$$\Delta^{s-r}b_+ =d_1\sigma_{i_1}d_2\sigma_{i_2}...d_k\sigma_{i_k}U \stackrel{k}{\rightarrow}\Delta^{s-r} \Delta^r =d_1d_2...d_kU, (E)$$

In $b_+\stackrel{k}{\rightarrow}\Delta^r$, each line of the sequence is noted
 by $K_i$: $b_+=K_k$, $\Delta^r=K_0$.If we note $\Delta^{s-r}Ki=L_i$,

$(E)$ as above gives $L_k=\Delta^{s-r}b_+ =d_1\sigma_{i_1}d_2\sigma_{i_2}...d_k\sigma_{i_k}U$, $L_0=\Delta^{s-r} \Delta^r =d_1d_2...d_kU$.

Finally, each $L_i$ is right divided by $\Delta^{s-r}$ and Lemma 11 1 gives,by supposing $s-r$ even,

$$\Delta^{s-r}b_+ =b_+\Delta^{s-r}=A_1\delta_{i_1}A_2\delta_{i_2}...A_k\delta_{i_k}A_{k+1}\Delta^{s-r} 
\stackrel{k}{\rightarrow}\Delta^r \Delta^{s-r} =A_1A_2...A_kA_{k+1} \Delta^{s-r}$$

By simplifying by $\Delta^{s-r}$, we have the result.

\vspace{0.2in}

We give now the proof of Lemma 3 1.

\vspace{0.1in}

PROOF Lemma 3 1:

Let $b=\Delta^r b_+$ be a quasipositive word and
$b_+=A_1\sigma_{i_1}A_2\sigma_{i_2}...A_k\sigma_{i_k}A_{k+1}\stackrel{k}{\rightarrow}\Delta^r =A_1A_2...A_kA_{k+1}$ be the corresponding k-sequence.

We have $b_+=b_1b_2...b_r...b_s$ in its left normal form.If $A_1$ as above is written $A_1=u_1...u_p$ in its left normal form, we have:

1-$A_1 \prec \Delta^r$ so $p \leq r \leq s$

2-$A_1 =u_1...u_p \prec b_1...b_p$

\vspace{0.1in}

We prove now that $A_1=b_1...b_p$.

By induction on $k$.

\vspace{0.1in}

If $k=1$. We deal with $b_+=A_1\sigma_{i_1}A_2$ and $\Delta^r=A_1A_2$.

In one side, we have $A_1\prec b_1...b_p$.

In the other side, Lemma 10 1 tells us we can consider $A_1$ equals to lgcd $(b_+,\Delta^r)$.
As, $p \leq r$, $b_1...b_p \prec \Delta^r$: so, $b_1...b_p \prec$ lgcd $(b_+,\Delta^r) =A_1 \prec b_1...b_p$.

Finally, $A_1=b_1...b_p$.

\vspace{0.1in}

For any $k$. We suppose the result true for $k-1$: for integers $r,s,p$,

with $p \leq r \leq s$, such as $b_1b_2...b_r..b_s=C_1\sigma_{i_1}C_2\sigma_{i_2}...C_{k-1}\sigma_{i_{k-1}}C_k \stackrel{k-1}{\rightarrow}\Delta^{r}=C_1C_2...C_k$ and $C_1=u_1u_2...u_p$ in its normal form , we have $C_1=b_1b_2...b_p$.

Let $b_+=b_1b_2...b_r...b_s=A_1\sigma_{i_1}A_2\sigma_{i_2}...A_k\sigma_{i_k}A_{k+1} \stackrel{k}{\rightarrow}\Delta^r=A_1A_2...A_kA_{k+1}$ be a k-sequence with
$A_1=u_1u_2...u_p \prec b_1b_2...b_p$.

Now, as $p \leq r$, $\Delta^{r}=b_1b_2...b_pUb_{t+1}...b_s$ for a positive $U$ and $s+p-r=t$.

The sequence $b_+=b_1b_2...b_pb_{p+1}...b_tb_{t+1}...b_s \stackrel{k}{\rightarrow}\Delta^{r}=b_1b_2...b_pUb_{t+1}...b_s$ is equivalent to

$$b_{p+1}...b_t \stackrel{k}{\rightarrow}U$$

we can write as $b_{p+1}...b_t =X_k\sigma_{i_k}Y_k\stackrel{1}{\rightarrow}X_kY_k\stackrel{k-1}{\rightarrow}U$.

By lemmas above, we obtain 

$b_1b_2...b_pb_{p+1}...b_tb_{t+1}...b_s =b_1b_2...b_pX_k\sigma_{i_k}Y_kb_{t+1}...b_s =A_1\sigma_{i_1}A_2\sigma_{i_2}...A_k\sigma_{i_k}A_{k+1}\stackrel{1}{\rightarrow}b_1b_2...b_pX_kY_kb_{t+1}...b_s=A_1\sigma_{i_1}A_2\sigma_{i_2}...A_kA_{k+1}\stackrel{k-1}{\rightarrow}b_1b_2...b_pUb_{t+1}...b_s$.

Finally the (k-1)-sequence $L_{k-1}=b_1b_2...b_pX_kY_kb_{t+1}...b_s=A_1\sigma_{i_1}A_2\sigma_{i_2}...A_kA_{k+1}\stackrel{k-1}{\rightarrow}L_0=b_1b_2...b_pUb_{t+1}...b_s$ with $A_1 \prec b_1b_2...b_p$ gives by induction

$$A_1=b_1b_2...b_p.$$

Moreover, in the previous sequence, we can notice that by Lemma 11 1

rgcd $(L_k,L_{k-1},...,L_0) \succ b_{t+1}...b_s $;as we can suppose $A_{k+1}=$ rgcd $(L_k,L_{k-1},...,L_0)$, $A_{k+1}=A_{k+1}' b_{t+1}...b_s$ for a positive $A_{k+1}'$.

We can now conclude by 

$b_{p+1}...b_t=\sigma_{i_1}A_2\sigma_{i_2}...A_k\sigma_{i_k}A_{k+1}'\stackrel{k}{\rightarrow}U=A_2A_3...A_kA_{k+1}'$.

\vspace{0.1in}

In the next paragraph, we give a computational algorithm and its complexity.

\paragraph{Algorithm}

We give now the practical consequence of Lemma 11 1.

\newtheorem{lemme12}{Lemma 12}
\begin{lemme12}

Let $b=\Delta^{-r}b_1b_2...b_r...b_s$ be a quasipositive word. By previous results, for $p$ such as $1 \leq p  \leq r \leq s$, there exist $k$ positive words 
$A_2,A_3,...,A_k,A_{k+1}$ and $k$ generators such as

$$b_{p+1}...b_{p+s-r}=\sigma_{i_1}A_2\sigma_{i_2}A_3...A_k\sigma_{i_k}A_{k+1}$$
$$U=A_2A_3...A_kA_{k+1}$$

where $U$ satisfies $b_1b_2...b_pUb_{p+s-r+1}...b_s=\Delta^r$.

Then, we can write 
$$b_{p+1}...b_{p+s-r}=\sigma_{i_1}'A_2'\sigma_{i_2}'A_3'...A_k'\sigma_{i_k}'A_{k+1}'$$
$$U=A_2'A_3'...A_k'A_{k+1}'$$

such as , for $2 \leq j \leq {k+1}$, $A_j' =$ lgcd $(A_j'\sigma_{i_j}'A_{j+1}'...A_k'\sigma_{i_k}'A_{k+1}',A_j'A_{j+1}'...A_k'A_{k+1}')$

\end{lemme12}

PROOF:

If for a k-sequence $\sigma_{i_1}A_2\sigma_{i_2}A_3...A_k\sigma_{i_k}A_{k+1}\stackrel{k}{\rightarrow}A_2A_3...A_kA_{k+1}$, there exist
$A_j'$ as in the lemma above, we say that the sequence satisfies (E).

We use here a double induction.Induction on $k$ and on $l=$ length $(A_2...A_{k+1})$.We say that E($k$,$l$) is true if  any k-sequence of
length $l$ satisfies (E).

The basis of induction are:

1-E($k$,$0$) true for any $k$.Indeed,$l=0$ means $A_j=1$ and the k-sequence is $\sigma_{i_1}\sigma_{i_2}...\sigma_{i_k}\stackrel{k}{\rightarrow}1$.
 Obviously, $A_j=1=$ lgcd $(\sigma_{i_1}\sigma_{i_2}...\sigma_{i_k},1)$.

2-E($1$,$l$) true for any $l$.The 1-sequence is $\sigma_{i_1}A_2\stackrel{1}{\rightarrow}A_2$ and obviously, $A_2=$ lgcd $(A_2,A_2)$.

We can suppose E($k$,$l-1$) true for any $k$ and E($k'$,$l$) true for $1 \leq k' \leq k-1$ and let us prove E($k$,$l$).

\vspace{0.1in}

Let $A_d$ be the first factor not equal to $1$, the k-sequence becomes 
$\sigma_{i_1}...\sigma_{i_d}A_d...A_k\sigma_{i_k}A_{k+1}\stackrel{k}{\rightarrow}A_d...A_{k+1}$.

As $A_d=\delta A_d'$, we get a (k+1)-sequence $X=\sigma_{i_1}...\sigma_{i_d}\delta A_d'...A_k\sigma_{i_k}A_{k+1}\stackrel{k}{\rightarrow}Y=A_d'...A_{k+1}$.

We have length $(A_d'...A_k)=l-1$ and as E($k$,$l-1$) true for ay $k$,

$$X=C_1\delta_{i_1}C_2\delta_{i_2}C_3...C_{k+1}\delta_{i_{k+1}}C_{k+2}$$

$$Y=C_1C_2...C_{k+1}C_{k+2}$$

Deleting $\delta$ corresponds to deleting 
$\delta_{i_d}$ so $\delta A_d'...A_{k+1}=\delta Y =\delta C_1C_2...C_{k+1}C_{k+2}=C_1...C_d\delta_{i_d}C_{d+1}...C_{k+2}$.

\vspace{0.2in}

1-If $C_1=\alpha C_1'$ is not equal to one,from the k-sequence of length $l$, 
$\sigma_{i_1}C_1\delta_{i_1}C_2\delta_{i_2}C_3...C_{k+1}\delta_{i_{k+1}}C_{k+2}\stackrel{k}{\rightarrow}C_1C_2...C_{k+1}C_{k+2}$, 
we extract the other one of length $(l-1)$ where all generators, but $\delta_{i_d}$,are
deleted

$$X'=\sigma_{i_1}C_1'\delta_{i_1}C_2\delta_{i_2}C_3...C_{k+1}\delta_{i_{k+1}}C_{k+2}\stackrel{k}{\rightarrow}Y'=C_1'C_2...C_{k+1}C_{k+2}.$$

E($k$,$l-1$) is true so $X'=D_1\delta_{i_1}'D_2\delta_{i_2}'D_3...D_k\delta_{i_k}'D_{k+1}$ and $Y'=D_1D_2...D_{k+1}$ with

$D_j =$ lgcd $(D_j\delta_{i_j}'D_{j+1}...D_k\delta_{i_k}'D_{k+1},D_jD_{j+1}...D_kD_{k+1})$.

Moreover, lgcd $(X,Y)=$ lgcd $(\alpha X',\alpha Y') = \alpha D_1$ and we obtain all datas for the final result, (F).

\vspace{0.1in}

2-If $C_1=1$.Then lgcd $(X,Y)=1$.We note $R=$ lgcd $(X,\delta Y)$ and let $\alpha $ such as $\alpha \prec R$.

\hspace{0.1in}2.1-if $\alpha $ different from $\delta $, then $\alpha \prec (X,\delta Y)$ gives $\alpha \prec (X,Y)$.As lgcd $(X,Y)=1$,
$\alpha = 1$.

\vspace{0.1in}

\hspace{0.1in}2.2-if $\alpha = \delta$, $\delta \prec X=\delta_{i_1}C_2\delta_{i_2}C_3...C_{k+1}\delta_{i_{k+1}}C_{k+2}$.

\vspace{0.1in}

\hspace{0.2in}2.2.1-if $\delta = \delta_{i_1}$ then $X=\delta C_2\delta_{i_2}C_3...C_{k+1}\delta_{i_{k+1}}C_{k+2}$ and
$Y\delta Y= \delta C_2...C_{k+1}C_{k+2}$.

As before, for $\sigma_{i_1} C_2\delta_{i_2}C_3...C_{k+1}\delta_{i_{k+1}}C_{k+2}\stackrel{k}{\rightarrow} C_2...C_{k+1}C_{k+2}$,

E($k$,$l-1$) beeing true, we have the final result, see (F).

\vspace{0.1in}

\hspace{0.2in}2.2.2-if $\delta$ is different from $\delta_{i_1}$ then $\delta \prec X=\delta_{i_1}C_2\delta_{i_2}C_3...C_{k+1}\delta_{i_{k+1}}C_{k+2}$

gives $\delta \prec X'=C_2\delta_{i_2}C_3...C_{k+1}\delta_{i_{k+1}}C_{k+2}$. Moreover, $\delta \prec \delta Y=C_2...C_d\delta_{i_d}C_{d+1}...C_{k+2}$.

As E($k-1$,$l$) is true, the (k-1)-sequence beeing $X'=C_2\delta_{i_2}C_3...C_{k+1}\delta_{i_{k+1}}C_{k+2}\stackrel{k-1}{\rightarrow}\delta Y$

where all $\delta_{i_j}$ but $\delta_{i_d}$ are deleted, we obtain

$$X'=\delta D_1\delta_{i_1}'D_2\delta_{i_2}'D_3...D_{k-1}\delta_{i_k}'D_k$$

$$\delta Y=\delta D_1...D_k$$

We use now Lemma 11 1.Given a k-sequence $L_k=A_1\sigma_{i_1}A_2\sigma_{i_2}A_3...A_k\sigma_{i_k}A_{k+1}\stackrel{k}{\rightarrow}L_0=A_1...A_{k+1}$

where $\sigma_{i_k}$ is first deleted,then $\sigma_{i_{k-1}}$...and so on,we can suppose rgcd $(L_k,...,L_0)=A_k$.

By symmetry, if we delete first $\sigma_{i_1}$, then $\sigma_{i_2}$...and so on, we can suppose $A_1=$ lgcd $(L_k,..,L_0)$.

Now, if we consider $L_k=\delta_{i_1} X'=\delta_{i_1}\delta D_1\delta_{i_1}'D_2\delta_{i_2}'D_3...D_{k-1}\delta_{i_k}'D_k
\stackrel{1}{\rightarrow}L_{k-1}=\delta D_1\delta_{i_1}'D_2\delta_{i_2}'D_3...D_{k-1}\delta_{i_k}'D_k
\stackrel{1}{\rightarrow}L_{k-2}=\delta D_1D_2\delta_{i_2}'D_3...D_{k-1}\delta_{i_k}'D_k
\stackrel{1}{\rightarrow}L_{k-3}=\delta D_1D_2D_3\delta_{i_3}'...D_{k-1}\delta_{i_k}'D_k
\stackrel{k-3}{\rightarrow}L_0=delta D_1D_2D_3...D_{k-1}D_k$, we can write $L_k$ as

$$L_k=\delta E_1\sigma_{i_1}''E_2\sigma_{i_2}''E_3...E_k\sigma_{i_k}''E_{k+1}$$
$$L_0=\delta E_1E_2E_3...E_kE_{k+1}$$

Finally, $X=L_k$ and $Y=L_0$ and doing as in (F), with $\sigma_{i_1} E_1\sigma_{i_1}''E_2\sigma_{i_2}''E_3...E_k\sigma_{i_k}''E_{k+1}
\stackrel{k}{\rightarrow}E_1E_2E_3...E_kE_{k+1}$, we obtain the final result.

\vspace{0.2in}

REMARK:

This last result gives practical consequences for the next algorithm.
Indeed, with the following datas,

$$P\stackrel{k}{\rightarrow}N$$
$$P=\sigma_{i_1}A_2\sigma_{i_2}A_3...A_k\sigma_{i_k}A_{k+1}$$
$$N=A_2A_3...A_kA_{k+1}$$

$\sigma_{i_1}$ is easy to compute, among the generators which begin $P$.

Moreover $A_2$ satisfies $A_2=$ lgcd $({\sigma_{i_1}}^{-1}P, N)$...and so on.


\vspace{0.2in}

\newtheorem{guess2}{Theorem 2}
\begin{guess2}

ALGORITHM

Let  $b$ be a word in $B_{n+1}$.For a positive word $w$, we note $beg(w)$ the set of generators left  dividing $w$,
$beg(w)=(\sigma / \sigma \prec w)$.

\vspace{0.1in}

1-We compute the algebraic length of $b$, $k$:

\hspace{0.1in}  if $k \leq 0$, STOP ($b$ is not quasipositive;if $k=0$, $b$ is quasipositive iff $b=1$)

\hspace{0.1in}  elsewhere we go to step 2

\vspace{0.1in}

2-We compute the left normal form of $b$, $b=\Delta^{-r}b_1...b_rb_{r+1}...b_s$:

\hspace{0.1in} if $r \leq 0$, STOP ($b$ is positive so quasipositive)
 
\hspace{0.1in} elsewhere we go to step 3

\vspace{0.1in}

3-For each $0 \leq j \leq r$, we compute $P_j=b_{j+1}...b_{s-r+j}$ and we compute $N_j$ such as $b_1b_2...b_jN_jb_{s-r+j+1}...b_s=\Delta^{r}$

\vspace{0.1in}

4-We consider the set S with S $=$ [$P_0,...,P_r,N_0,...,N_r$] and the empty set S'

\vspace{0.1in}

5-We begin with $j=0$

\vspace{0.1in}

6-We set $P=P_j$ and $N=N_j$ and we add $P$ to the set S'

\vspace{0.1in}

7-We set  $m=1$

\vspace{0.1in}

8-For all elements $P$ in the set S':

 we compute $beg(P)$ and for each $\sigma$ in $beg(P)$, we compute $P'$ such as $P=\sigma P'$.

For each $P'$, we compute the normal form of $N^{-1}P'$, we can write as $N^{-1}P'=A^{-1}B$.

We add all theses $P',A,B$ in the set $S$ and we add all the $B$ in the set S'.

\vspace{0.1in}

10-we set $m:=m+1$

\vspace{0.1in}

9-If $m=k$, we set $j:=j+1$:

\hspace{0.1in} if $j=r+1$ then we go to step 10

\hspace{0.1in} elsewhere, we go to step 6

Elsewhere, if $m\leq k-1$, we go to step 8.

\vspace{0.1in}

10-If $1$ belongs to the set S, then $b$ is quasipositive.Keeping
trace to the factors used,by noting them by $A_2,...,A_{k+1}$ and $\sigma_{i_1},...,\sigma_{i_k}$, we output the quasipositive form of 
a conjuguate of $b$:

$$b'=A_{k+1}^{-1}...A_2^{-1}\sigma_{i_1}A_2\sigma_{i_3}...\sigma_{i_k}A_{k+1}$$

and the form of $b$, $b=T_j^{-1}b'T_j$ with $T_j=(b_{s-r+j+1}...b_s)$.

Elsewhere, STOP ($b$ is not quasipositive).

\end{guess2}

EXAMPLES:

In $B_4$, we put $b=2212^{-1}2^{-1}2^{-1}2^{-1}322$ with $k=2$;its normal form is

$$b = \Delta^{-4} 2321.32123.12321.32123.12.2132.2$$

We have $r=4$, $s=6$ and $b_1=2321$, $b_2=32123$, $b_3=13221$, $b_4=12$, $b_5=2132$, $b_6=2$

so $a_1=23$, $a_2=2$, $a_3=2$, $a_4=1321$, $a_5=13$, $a_6=13213$.

We have $P=b_3b_4b_5 = 12321.32123.12$;the obtained $N$ is $N=a_2\bar{a_1}a_6\bar{a_5} = 2.21.13213.31=2.21.132113.3$.

But $P=123231212312=122321212312=122312112312=122132112132=122132121232=1.2213.1211232=1.221132113.2.3$.

We have $P=1.A_2.2.A_3$ with $A_2=221132113$ and $A_3=3$: as $N=A_2A_3$ $b$ is quasipositive.

\vspace{0.1in}

If we consider now $b=123^{-1}$ with $k=1$;its normal form is 

$b= \Delta^{-1} 323.2132$ so $b_1=323$,$b_2=2132$ and $r=1$; and $a_1=123$, $a_2=13$.

We have two possibilities for $P$:

$P=b_2 =2132$ and $N=a_1 = 123$: $P$ can't be reduced to $N$

$P=b_1=323$ and $N=\bar{a_2}=31$: $P$ can't be reduced to $N$

Finally $b$ is not quasipositive.

\vspace{0.2in}

If we consider now $b$, a word of algebraic length $k$ written in its left normal form $b=\Delta^{-r}b_+$, and we note by $l$ the length of $b_+$,we obtain the next result.

\newtheorem{guess3}{Theorem 3}
\begin{guess3}
The complexity of the previous algorithm is $O(n^kl^2)$.

\end{guess3}

PROOF:

For the initial $P_j$ and $N_j$, we have $r+1$ choices.
Let us note by $x$ the complexity for computing the normal form of a word.

At the end of the algorithm, the set S contains at most $(r+1)n^k$ elements: for an initial $P_j$, we have to compute $beg(P)$  $k$ times and we have $n$ 
choices of generators for this last.

Moreover at each step, we have to compute the normal form of the elements we are dealing with(step 8), so the final complexity is about 

$(r+1)n^kx$.

In different litterature [3], it's proved that $x$ is about $nl^2$.

Finally, the final complexity is at most $n^{k+1}l^2$, which is  $O(n^kl^2)$.

\paragraph{Conclusion}

This algorithm and all  its proofs  only need  the combinatorial properties of  braid groups: the existence 
of lgcd, rgcd and $\Delta$, and as a consequence,the existence of a normal form for every word and over all, the 
existence of a length function.
Theses conditions beeing the classical definition of  Garside groups,as  braid groups,spherical Artin-Tits groups, (see [2],[3],[4]), 
the previous algorithm can then  be used in all theses groups,with a same complexity.

\paragraph{References}

[1]BIRMAN, Joan S. Braids, links, and mapping class groups. {\it Annals of Mathematics Studies}, No. 82.
 Princeton University Press, Princeton, N.J.; University of Tokyo Press, Tokyo, 1974.

[2]BRIESKORN,SAITO.Artin gruppen und Coxeter gruppen.{\it Invent. Math. 17}, 1972, 245-271

[3]DEHORNOY.Groupes de Garside.{\it Ann. Sc. Ec. Norm. Sup.}, No. 35, 2002,267-306

[4]DEHORNOY,PARIS.Gaussian groups and Garside groups, two generalisatios of Artin groups.{\it Proc. London. Math. Soc(3)}, No 79, 1999, 569-604.

[5]GARSIDE, F. A. The braid group and other groups. {\it Quart. J. Math. Oxford Ser+. (2)  20  1969 235--254.

[6]EL-RIFAI, ELSAYED A.; MORTON, H. R. Algorithms for positive braids. {\it Quart. J. Math. Oxford Ser}, No 45  (1994),  no. 180, 479--497.

[7]OREVOV, S. Yu. Erratum to: "Link theory and oval arrangements of real algebraic curves" [Topology 38 (1999), no. 4, 779--810;
 MR 2000b:14066]. {\it Topology} 41 (2002), no. 1, 211--212. 14P25 (20F36 57M25). 34 (2000), no. 3, 84--87.

[8]OREVKOV, S. Yu. Quasipositivity test via unitary representations of braid groups and its applications to real algebraic curves.
 {\it J. Knot Theory Ramifications} 10 (2001), no. 7, 1005--1023.

[9]RUDOLPH, Lee Quasipositive pretzels.  {\it Topology Appl.}  115  (2001),  no. 1, 115--123.

[10]RUDOLPH, Lee Quasipositive plumbing (constructions of quasipositive knots and links. V). {\it Proc. Amer. Math. Soc.} 126 (1998),
 no. 1, 257--267.

\vspace{0.2in}

\it{A.Bentalha}

\it{Universite Paul abatier}

\it{Laboratoire Emile Picard, 31400 Toulouse, France}

\it{a.bentalha@laposte.net}

\end{document}